\documentclass{pnastwo}

\usepackage{graphicx}

\usepackage{amssymb,amsfonts,amsmath,bbm}

\newtheorem{thm}{Theorem}

\let\set\mathbbm
\def\<#1>{\langle#1\rangle}

\def\qbinom#1#2{\genfrac{[}{]}{0pt}{}{#1}{#2}}

\contributor{Edited by George E. Andrews, Pennsylvania State University, University Park, PA,
and approved December 22, 2010 (received for review February 24, 2010)}
\url{www.pnas.org/cgi/doi/10.1073/pnas.1019186108}
\copyrightyear{2008}
\issuedate{February 8, 2011}
\volume{vol. 108}
\issuenumber{no. 6}
\footlineauthor{Koutschan, Kauers, Zeilberger}
\begin{document}

%%%%%%%%%%%%%%%%%%%%%%%%%%%%%%

\title{Proof of George~Andrews's and David~Robbins's q-TSPP~conjecture}

%% Enter authors via the \author command.  
%% Use \affil to define affiliations.
%% (Leave no spaces between author name and \affil command)

%% Note that the \thanks{} command has been disabled in favor of
%% a generic, reserved space for PNAS publication footnotes.

%% \author{<author name>
%% \affil{<number>}{<Institution>}} One number for each institution.
%% The same number should be used for authors that
%% are affiliated with the same institution, after the first time
%% only the number is needed, ie, \affil{number}{text}, \affil{number}{}
%% Then, before last author ...
%% \and
%% \author{<author name>
%% \affil{<number>}{}}

%% For example, assuming Garcia and Sonnery are both affiliated with
%% Universidad de Murcia:
%% \author{Roberta Graff\affil{1}{University of Cambridge, Cambridge,
%% United Kingdom},
%% Javier de Ruiz Garcia\affil{2}{Universidad de Murcia, Bioquimica y Biologia
%% Molecular, Murcia, Spain}, \and Franklin Sonnery\affil{2}{}}

\author{%
  Christoph Koutschan\affil{1}{Department of Mathematics,
    Tulane University, New Orleans, LA 70118}\footnotemark\ ,
  Manuel Kauers\affil{2}{Research Institute for Symbolic Computation,
    Johannes Kepler University Linz, A 4040 Linz, Austria}\footnotemark\ ,
  Doron Zeilberger\affil{3}{Mathematics Department, Rutgers University (New Brunswick), Piscataway, NJ 08854-8019}
}

%\contributor{Submitted to Proceedings of the National Academy of Sciences
%of the United States of America}

%% The \maketitle command is necessary to build the title page.
\maketitle

\def\thefootnote{\fnsymbol{footnote}}

\footnotetext[2]{Present address: Research Institute for Symbolic Computation, Johannes Kepler University Linz, A 4040 Linz, Austria.}
\footnotetext[4]{To whom correspondence should be addressed. E-mail: mkauers@risc.jku.at.}

%%%%%%%%%%%%%%%%%%%%%%%%%%%%%%%%%%%%%%%%%%%%%%%%%%%%%%%%%%%%%%%%
\begin{article}

\begin{abstract} 
The conjecture that the orbit-counting generating function for totally
symmetric plane partitions can be written as an explicit
product formula, has been stated independently by George Andrews and
David Robbins around 1983.  We present a proof of this long-standing
conjecture.
\end{abstract}

%% When adding keywords, separate each term with a straight line: |
\keywords{computer algebra | enumerative combinatorics | partition theory}

%% Optional for entering abbreviations, separate the abbreviation from
%% its definition with a comma, separate each pair with a semicolon:
%% for example:
%% \abbreviations{SAM, self-assembled monolayer; OTS,
%% octadecyltrichlorosilane}

% \abbreviations{}

%% The first letter of the article should be drop cap: \dropcap{}
%\dropcap{I}n this article we study the evolution of ''almost-sharp'' fronts

%% Enter the text of your article beginning here and ending before
%% \begin{acknowledgements}
%% Section head commands for your reference:
%% \section{}
%% \subsection{}
%% \subsubsection{}

\section{1\kern.5em Proemium}

%\dropcap 
In the historical conference {\it Combinatoire \'{E}numerative} that took place
at the end of May 1985 in Montr\' eal,
Richard Stanley raised some intriguing problems about the
enumeration of plane partitions (see below), which he later expanded into a
fascinating article~\cite{Sta1}. Most of these problems concerned the
enumeration of ``symmetry classes'' of plane partitions that were
discussed in more detail in another article of Stanley~\cite{Sta2}.
All of the conjectures in the latter article have since been proved
(see David Bressoud's modern classic~\cite{B}), except one, which until now 
has resisted the efforts of some of the greatest minds in enumerative
combinatorics.  It concerns the proof of an explicit formula for the
$q$-enumeration of totally symmetric plane partitions, conjectured
around 1983 independently by George Andrews
and David Robbins~(\cite{Sta2}, \cite{Sta1} conj.~7, \cite{B}
conj.~13, and already alluded to in \cite{Andrews80}).  In the 
present article we finally turn this conjecture into a theorem.

A plane partition $\pi$ is an array $\pi=(\pi_{i,j})_{1\leq i,j}$, of
non-negative integers $\pi_{i,j}$ with finite sum $|\pi|=\sum \pi_{i,j}$,
which is weakly decreasing in rows and columns so that
$\pi_{i,j}\geq\pi_{i+1,j}$ and $\pi_{i,j}\geq\pi_{i,j+1}$. A plane
partition~$\pi$ is identified with its 3D Ferrers diagram which is
obtained by stacking $\pi_{i,j}$ unit cubes on top of the
location~$(i,j)$. The result is a \hbox{left-,} \hbox{back-,} and bottom-justified
structure in which we can refer to the locations $(i,j,k)$ of the
individual unit cubes. If the diagram is invariant under the action of
the symmetric group~$S_3$ on the coordinate axes then $\pi$ is called a totally symmetric
plane partition (TSPP). In other words, $\pi$ is called totally
symmetric if whenever a location $(i,j,k)$ in the diagram is occupied
then all its up to~$5$ permutations
$\{(i,k,j),(j,i,k),(j,k,i),(k,i,j),(k,j,i)\}$ are occupied as well.
Such a set of cubes, i.e., all cubes to which a certain cube can be
moved via $S_3$ is called an orbit; the set of all orbits of~$\pi$
forms a partition of its diagram (see Figure~1). %\ref{tsppfigure}).

In 1995, John Stembridge~\cite{Ste} proved Ian Macdonald's conjecture
that the number of totally symmetric plane partitions with largest part at most~$n$, i.e.,
those whose 3D Ferrers diagram is contained in the cube $[0,n]^3$,
is given by the elegant product-formula
\[
  \prod_{ 1 \leq i \leq j \leq k \leq n} \frac{i+j+k-1}{i+j+k-2}.
\]
Ten years after Stembridge's completely human-generated proof,
Andrews, Peter Paule and Carsten Schneider~\cite{APS} came up with a
computer-assisted proof based on an ingenious matrix decomposition,
but since no $q$-analog of their decomposition was found, their
proof could not be extended to a proof for the $q$-case.  A third
proof of Stembridge's theorem~\cite{Koutschan09,Koutschan10a}, even more
computerized, was recently found in the context of our
investigations of the $q$-TSPP conjecture. 
We have now succeeded in
completing all the required computations for an analogous proof of the
$q$-TSPP conjecture, and can therefore announce:

\begin{figure}
\centerline{\includegraphics[width=0.5\textwidth]{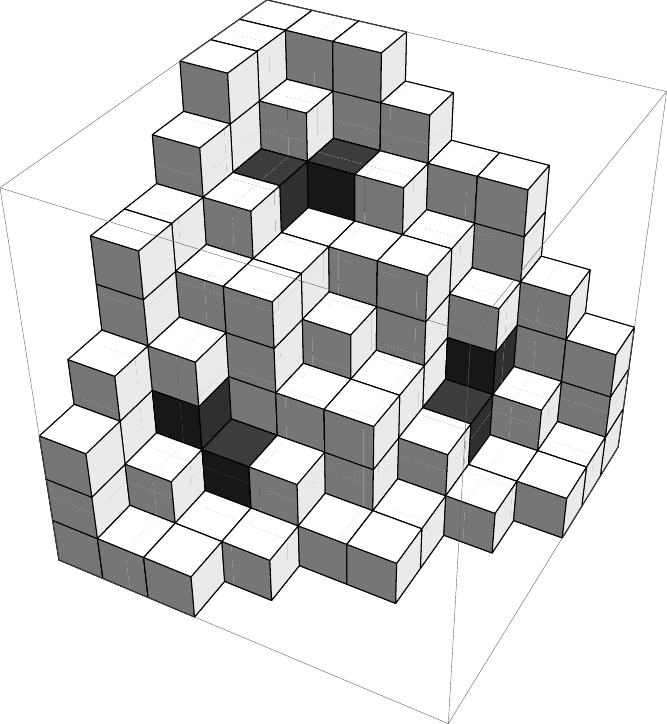}}
\caption{A totally symmetric plane partition with largest
part eight. The dark cubes form an orbit of size six (corresponding
to all permutations of (2,3,6)). Each cube belongs to exactly 
one orbit. An orbit may consist of one, three, or six cubes.}\label{tsppfigure}
\end{figure} 

\begin{thm}\label{qTSPPthm}
Let $\pi/S_3$ denote the set of orbits of a totally symmetric plane partition~$\pi$
under the action of the symmetric group~$S_3$.  Then the
orbit-counting generating function (\cite[p.~200]{B},
\cite[p.~106]{Sta2}) is given by
\[
  \sum_{\pi\in T(n)} q^{|\pi/{S_3}|}=\prod_{1 \leq i \leq j \leq k \leq n}\frac{1-q^{i+j+k-1}}{1-q^{i+j+k-2}}
\] 
where $T(n)$ denotes the set of totally symmetric plane partitions
with largest part at most~$n$.
\end{thm}

\smallskip
\noindent\emph{Proof sketch.}
 Our proof is based on a result by Soichi Okada~\cite{O} who has shown that
 the theorem is implied by a certain---conjectured---determinant
 evaluation.  These preliminaries are stated accurately in the next section,
 followed by a description of the holonomic ansatz~\cite{Z3} that we
 are going to pursue.  This approach relies on a kind of oracle that
 tells us a description of a certificate function~$c_{n,j}$;
 the odyssey how this function has been ``guessed'' is described in the following section.
 Once $c_{n,j}$ is known, the determinant evaluation reduces further to
 proving the three identities \eqref{eq:Normalization},
 \eqref{eq:Soichi}, and \eqref{eq:Okada} stated below.  
 Finally we come full circle by proving these identities.
 Technical details of the proof and our computations 
 as well as the explicit certificates are provided electronically 
 on our Web site http:/$\!$/www.risc.jku.at/people/ckoutsch/qtspp/.
 Further details which are suppressed here can be found in the detailed
 description of the proof for the case $q=1$ which proceeds along the same 
 lines~\cite{Koutschan09,Koutschan10a}. \rule{1ex}{1ex}
\par\medskip

Our proof is noteworthy not only for its obvious significance in enumerative combinatorics, where
it settles a long-standing conjecture, attempted by many people. It is noteworthy also for
computational reasons, as the computations we performed went far beyond what has
been thought to be possible with currently known algebraic algorithms,
software packages, and computer hardware.

\section{2\kern.5em The Telemachiad}\label{secReduction}

In order to prove the $q$-TSPP conjecture, we exploit an elegant reduction
by Okada~\cite{O} to the problem of evaluating a certain determinant.
This determinant is also listed as Conjecture~46 in Christian Krattenthaler's
essay~\cite{K} on the art of determinant evaluation.

Let, as usual, $\delta_{i,j}$ be the Kronecker delta function
and let, also as usual,
\[
  \qbinom n k = \frac{(1-q^n)(1-q^{n-1}) \cdots (1-q^{n-k+1})}{(1-q^k)(1-q^{k-1}) \cdots (1-q)}
\]
denote the $q$-binomial coefficient.
Define the discrete function $a_{i,j}$ for $i,j\geq1$ by
\[%\textstyle
 % a_{i,j}:=  
  q^{i+j-1} \left(\qbinom{i+j-2}{i-1} + q \qbinom{i+j-1}{i} \right)+
  (1+q^i)\delta_{i,j}-\delta_{i,j+1}.
\]
Okada's crucial insight is that Theorem~\ref{qTSPPthm} holds if
\[
  \det(a_{i,j})_{1 \leq i,j \leq n}=
  \!\!\!\!\prod_{ 1 \leq i \leq j \leq k \leq n}\!\! {\left ( \frac{ 1-q^{i+j+k-1}}{1-q^{i+j+k-2}}\right )^2}=:b_n
  \quad(n\geq 1).
\]
So for proving the $q$-TSPP conjecture, it is sufficient to prove this
conjectured determinant evaluation. For this purpose we apply a
computational approach originally proposed in~\cite{Z3} which is 
applicable to identities of the form
\[
  \det (a_{i,j})_{1 \leq i,j \leq n} = b_n\qquad(n\geq 1)
\] 
where $a_{i,j}$ and $b_n$ (with $b_n\neq0$ for all $n\geq 1$) are given
explicitly (as it is the case here).

The approach rests on the following induction argument on~$n$.  For
$n=1$, the identity is trivial.  Suppose the identity holds for $n-1$;
then the linear system
\[
  \begin{pmatrix}
    a_{1,1} & \cdots & a_{1,n-1} & a_{1,n} \\
    \vdots & \ddots & \vdots & \vdots \\
    a_{n-1,1} & \cdots & a_{n-1,n-1} & a_{n-1,n} \\
    0 & \cdots & 0 & 1
  \end{pmatrix}
  \begin{pmatrix}
    c_{n,1} \\
    \vdots \\
    c_{n,n-1} \\
    c_{n,n}
  \end{pmatrix}
  = 
  \begin{pmatrix}
    0 \\
    \vdots\\
    0\\
    1
  \end{pmatrix}
\]
has a unique solution $(c_{n,1},\dots,c_{n,n})$.
The component $c_{n,j}$ of this solution 
is precisely the $(n,j)$-cofactor divided by the $(n,n)$-cofactor 
of the $n\times n$ determinant. The division is meaningful because
the $(n,n)$-cofactor is just the $(n-1)\times(n-1)$ determinant,
which by induction hypothesis is equal to~$b_{n-1}$, which by general
assumption is nonzero. 
Since the $n\times n$ determinant can be expressed in terms of
the matrix entries $a_{n,j}$ and the normalized cofactors $c_{n,j}$
via
\[
  b_{n-1}\sum_{j=1}^n c_{n,j} a_{n,j},
\]
the induction step is completed by showing that this sum evaluates to~$b_n$.

The difficulty is that this last summation involves the
function~$c_{n,j}$ (of the discrete variables $n$ and $j$) for which we do not have an explicit expression for
general~$n$ and~$j$.  In order to achieve our goal, we guess a
suitable description of a function $c_{n,j}$ and then prove that it
satisfies the three identities
\begin{alignat}3
  c_{n,n}&=1&\qquad&(n\geq1),\label{eq:Normalization}\\
  \sum_{j=1}^n c_{n,j}a_{i,j}&=0&&(1\leq i<n),\label{eq:Soichi}\\
  \sum_{j=1}^n c_{n,j}a_{n,j}&=\frac{b_n}{b_{n-1}}&&(n\geq1)\label{eq:Okada}.
\end{alignat}
Once these identities are proved, then, by the argument given before the~$c_{n,j}$ must be 
precisely the normalized $(n,j)$-cofactors of the $n\times n$ determinant and
the determinant evaluation follows as a consequence. 
So in a sense, the function~$c_{n,j}$ plays the r\^{o}le of a
certificate for the determinant identity.

\section{3\kern.5em The Odyssey}\label{secGuess}

In our setting the certificate function $c_{n,j}$ will be described
implicitly by a system of linear recurrence equations in $n$ and $j$ with
coefficients depending polynomially on $q$, $q^j$ and~$q^n$.  Such recurrence
equations can be phrased as the elements of some noncommutative operator
algebras such as $\set Q(q,q^n,q^j)[S_n,S_j]$ where the symbols $S_x$ represent the
shifts $x\mapsto x+1$.  If a function is annihilated by certain operators (viz.\
it satisfies certain recurrence equations), then it is also annihilated by all
the elements in the (left) ideal generated by those operators.  We speak of an
annihilating ideal and represent such ideals by (left) Gr\"obner bases, so
that for instance ideal membership can be decided effectively.

The annihilating ideals we use for representing functions are such that they
uniquely determine the function up to some finitely many initial values which
we can list explicitly.
Technically, this requirement means that the ideals have dimension zero and that some
particular polynomial coefficients appearing in the recurrence system must not
vanish simultaneously for infinitely many points~$(n,j)$.  These recurrence
systems are similar but somewhat simpler than $q$-holonomic systems~\cite{Z1},
which satisfy some additional requirements that are not needed for our proof.
We have checked that all the ideals arising in our proof are indeed of the
desired form, but in the interest of clarity we suppress a more detailed
description of these checks here. For a complete analysis including all
technical details we refer to the supplementary material on our Web site.  Also
in the interest of clarity, we will from now on identify recurrence equations
with their corresponding operators. The details are to a large extent analogous to the special case $q=1$, which as mentioned in the introduction has been written up in detail in~\cite{Koutschan09,Koutschan10a}.

A priori, there is no reason why the normalized cofactors $c_{n,j}$ should admit
a recursive description of the kind we are aiming at, but there is also no reason why they
should not. It turns out that they do, and this situation is fortunate because
for functions described in this way, techniques are known by which the
required identities~\eqref{eq:Normalization}, \eqref{eq:Soichi} and
\eqref{eq:Okada} can be proven algorithmically~\cite{Z1,T,CS,C}.
In order to find a recursive description for~$c_{n,j}$, we first
computed explicitly the normalized cofactors $c_{n,j}\in\set Q(q)$ for a few hundred
specific indices $n$ and $j$ by directly solving the linear system
quoted in the previous section.  Using an algorithm reminiscent
of polynomial interpolation, we then constructed a set of recurrences
compatible with the values of $c_{n,j}$ at the (finitely many) indices
we computed.

Polynomial interpolation applied to a finite sample $u_1,\dots,u_k$ of
an infinite sequence $u_n$ will always deliver some polynomial~$p$ of
degree at most $k-1$ which matches the given data. If it turns out
that this polynomial matches some further sequence terms
$u_{k+1},u_{k+2},\dots$, then it is tempting to conjecture that
$u_n=p(n)$ for all~$n$.  The more specific points~$n$ are found to
match, the higher is the evidence in favor of this conjecture.

Very much analogously, it is possible to extract recurrence equations
from some finite number of values~$c_{n,j}$. The equations become
trustworthy if they also hold for points~$(n,j)$ which were not used for
their construction.  In this way, we have discovered a system of
potential recurrence equations for the~$c_{n,j}$, which, despite being
respectable in size (about 30 Mbytes), appears to be a rather
plausible candidate for a recursive description of the normalized
cofactors~$c_{n,j}$. The system is available for download on our
Web site. It consists of three equations involving shifts of order up to 
four (with respect to both $n$ and~$j$), with polynomial coefficients 
of degrees up to 52 (with respect to~$q^n$) and 24 (with 
respect to~$q^j$). Some further statistics on the set of operators 
can be found in a preliminary version of this article~\cite{Kauers}.

\section{4\kern.5em The Nostos}\label{secProof}

Now we switch our point of view. We discard the definition that
$c_{n,j}$ be the normalized $(n,j)$-cofactor and redefine the discrete
function $c_{n,j}$ as the unique solution of the guessed recurrence
system whose (finitely many) initial values agree with the normalized
$(n,j)$-cofactor. If we succeed in proving that $c_{n,j}$ defined in
this way satisfies the identities~\eqref{eq:Normalization},
\eqref{eq:Soichi} and \eqref{eq:Okada}, then we are done.
For this purpose we
provide operators which belong to the annihilating
ideal of $c_{n,j}$ or related ideals and have certain features which imply
the desired identities. So in a sense, these operators play the r\^ole
of certificates for the identities under consideration. 

Because of their astronomical size (up to 7~Gbytes; equivalent to
more than one million printed pages; corresponding to about 2.5~tons
of paper), these certificates are not included explicitly in this
article but provided only electronically on our Web site. Also because
of their size, it was not possible to construct them by
simply applying the standard algorithms from~\cite{Z1,T,C,CS}. A
detailed explanation of how exactly we found the certificates is
beyond the scope of this article and will be given in a separate
publication.  But this lack of explanation does not at all affect the
soundness of our proof,
because the correctness of the certificates can be checked
independently by simply performing ideal membership tests.  
Our certificates are so big that even this ``simple'' calculation is
not quite trivial, but a reader with a sufficient amount of patience
and programming expertise will be able to do it.

We proceed by explaining the properties of the certificates provided
on our Web site and why they imply \eqref{eq:Normalization},
\eqref{eq:Soichi} and~\eqref{eq:Okada}. 

\subsection{A certificate for {\normalsize\eqref{eq:Normalization}}}

To certify that $c_{n,n}=1$ for all $n\geq 1$, we provide a recurrence in
the annihilating ideal of $c_{n,j}$ which is of the special form
\begin{alignat*}1
  &p_7(q,q^j,q^n)c_{n+7,j+7}+p_6(q,q^j,q^n)c_{n+6,j+6}+\cdots\\
  &\quad\cdots+p_1(q,q^j,q^n)c_{n+1,j+1}+p_0(q,q^j,q^n)c_{n,j}=0.
\end{alignat*}
By virtue of the substitution $j\mapsto n$, it translates into a recurrence for
the diagonal sequence~$c_{n,n}$.  It is a relatively cheap computation to check
that this recurrence contains the operator $S_n-1$ annihilating the constant
sequence~$1$ as a (right) factor. This fact implies that $c_{n,n}$ and the constant
sequence~$1$ both satisfy the same seventh-order recurrence. Therefore, after
checking $c_{1,1}=c_{2,2}=\cdots=c_{7,7}=1$ and observing that 
$p_7(q,q^n,q^n)\neq0$ for all $n\in\set N$, it can be concluded that $c_{n,n}=1$ 
for all~$n$.

Similarly, we showed that $c_{n,0}=0$ for all $n\geq 1$ and that
$c_{n,j}=0$ for all $j>n$. This knowledge greatly simplifies the
following proofs of the summation identities.

\subsection{Certificates for {\normalsize\eqref{eq:Soichi}}} 

In order to prove the first summation identity, we
translate~\eqref{eq:Soichi} into the equivalent formulation
\begin{alignat*}1
  &\sum_{j=1}^n \frac{q^{i+j-1}(q^{i+j}+q^i-q-1)}{q^i-1}\qbinom{i+j-2}{i-1}c_{n,j}\\
  &=c_{n,i-1}-(q^i+1)c_{n,i}\qquad(1\leq i<n),
  \tag{$\mathbf{\ref{eq:Soichi}'}$}\label{Soichi1}
\end{alignat*}
taking into account that $c_{n,0}=0$.  We provide certificates
for~\eqref{Soichi1}.

To this end, let~$c'_{n,i,j}$ be the summand of the sum on the
left-hand side.  A recursive description for $c'_{n,i,j}$ can be
computed directly from the defining equations of $c_{n,j}$ and the
fact that the rest of the summand is a $q$-hypergeometric factor.  In
the corresponding operator ideal~$I'$ we were able to find two
different recurrence equations for $c'_{n,i,j}$ which are of the special
form
\begin{alignat*}1
  &p_{0,3}c'_{n,i+3,j}+
   p_{1,2}c'_{n+1,i+2,j}+
   p_{0,2}c'_{n,i+2,j}\\
  &+p_{2,0}c'_{n+2,i,j}+
   p_{0,1}c'_{n,i+1,j}+
   p_{0,0}c'_{n,i,j}\\
  &=u_{n,i,j+1}-u_{n,i,j}
\intertext{and}
  &r_{4,0}c'_{n+4,i,j}+
   r_{2,1}c'_{n+2,i+1,j}+
   r_{0,2}c'_{n,i+2,j}\\
  &+r_{2,0}c'_{n+2,i,j}+
   r_{0,1}c'_{n,i+1,j}+
   r_{0,0}c'_{n,i,j}\\
  &=v_{n,i,j+1}-v_{n,i,j}
\end{alignat*}
respectively, where the $p_{\mu,\nu}$ and $r_{\mu,\nu}$ are certain rational functions in
$\set{Q}(q,q^i,q^n)$ and $u_{n,i,j}$ and $v_{n,i,j}$ are
$\set{Q}(q,q^i,q^j,q^n)$-linear combinations of certain shifts of
$c'_{n,i,j}$ which are determined by the Gr\"obner basis of~$I'$
as described in~\cite{C}. Next observe that $c'_{n,i,j}=0$ for $j\leq 0$ (because
of the $q$-binomial coefficient) and also for $j>n$ (because
$c_{n,j}=0$ for $j>n$). When summing the two recurrence equations for $j$ from $-\infty$
to $+\infty$, the right-hand side telescopes to~$0$ and the left-hand
side turns into a recurrence for the sum in~\eqref{Soichi1}. This
method is known as creative telescoping.  Finally we obtain two
annihilating operators~$P_1$ and $P_2$ for the left-hand side
of~\eqref{Soichi1}.

For the right-hand side of \eqref{Soichi1}, we can again construct an ideal of
recurrences from the defining equations of~$c_{n,j}$.  It turns out that this
ideal contains $P_1$ and~$P_2$, so that both sides of \eqref{Soichi1} are
annihilated by these two operators. Additionally, $P_1$ and $P_2$ have been
constructed such that they require only finitely many initial values to produce
a uniquely determined bivariate sequence (and this finite number of needed
values is known explicitly). The proof was completed by checking that the two
sides of \eqref{Soichi1} agree at these finitely many points~$(n,j)$.

\subsection{The certificate for {\normalsize\eqref{eq:Okada}}}

To certify the final identity, we rewrite \eqref{eq:Okada} equivalently into
\begin{alignat*}1
 &(1+q^n)-c_{n,n-1}\\
 &\quad+\sum_{j=1}^n \frac{q^{n+j-1}(q^{n+j}+q^n-q-1)}{q^n-1}\qbinom{n+j-2}{n-1}c_{n,j}\\
 &\qquad=\frac{(q^{2n};q)_n^2}{(q^n;q^2)_n^2}
 \qquad(n\geq 1)\tag{$\mathbf{\ref{eq:Okada}'}$}\label{Okada1}
\end{alignat*}
where $(a;q)_n:=(1-a)(1-qa)\cdots(1-q^{n-1}a)$ denotes the
$q$-Pochhammer symbol.  As before, we use creative telescoping to
provide a certified operator~$P$ which annihilates the sum on the
left-hand side. In the present case, a single operator is sufficient
because the sum depends only on a single variable~$n$ (there is no~$i$
there).  Using the operator~$P$ and the defining equations for~$c_{n,j}$, 
we then construct a recurrence for the entire left-hand
side, which turns out to have order~$12$. This recurrence is a left multiple
of the second order operator 
\begin{alignat*}1
  &(q^n+1)^2 (q^{n+1}+1)^2 (q^{2n+1}-1)^2 (q^{2n+3}-1)^2 S_n^2\\
  &\quad-(q^{n+1}+q^{2n+2}+1)^2 (q^{3 n+1}-1)^2 (q^{3n+5}-1)^2 
\end{alignat*}
which can be seen to annihilate the right-hand side, and its leading coefficient does not vanish for any index~$n$.
Therefore the proof of Theorem~\ref{qTSPPthm} is completed by
checking~\eqref{Okada1} for $n=1,2,\dots,12$. \emph{Quod erat demonstrandum.}

\section{5\kern.5em Epilogue}\label{Conclusion}

Paul Erd{\H o}s famously believed that every short and elegant mathematical statement
has a short and elegant ``proof from The Book'', and if humans tried hard enough, they
would eventually find it. Kurt G\"odel, on the other hand, meta-proved that there exist
many short and elegant statements whose shortest possible proof is very long. It is very possible that
the $q$-TSPP theorem {\it does} have a yet-to-be-found {\it proof from The Book}, but it
is just as possible that it does not, and while we are sure that the present proof is not
the shortest possible (there were lots of random choices in designing the proof), it may well
be the case that the shortest-possible proof is still very long, and would still require heavy-duty
computer calculations.

% Be that as it may, we believe that our general approach, and our way of taming the
% computer to prove something that seemed intractable with today's hardware and software are
% very elegant, and deserves to be included in The Book.

%% == end of paper:

%% Optional Materials and Methods Section
%% The Materials and Methods section header will be added automatically.

%% Enter any subheads and the Materials and Methods text below.
%\begin{materials}
% Materials text
%\end{materials}

%% Optional Appendix or Appendices
%% \appendix Appendix text...
%% or, for appendix with title, use square brackets:
%% \appendix[Appendix Title]

\begin{acknowledgments}
C.K. was supported by grants National Science Foundation---Division of Mathematical Sciences 0070567
 and Fonds zur wissenschaftlichen F\"orderung (FWF) P20162-N18.
M.K. was supported by FWF Y464-N18.
D.Z. was supported in part by the US National Science Foundation.
\end{acknowledgments}

%% PNAS does not support submission of supporting .tex files such as BibTeX.
%% Instead all references must be included in the article .tex document. 
%% If you currently use BibTeX, your bibliography is formed because the 
%% command \verb+\bibliography{}+ brings the <filename>.bbl file into your

\begin{thebibliography}{10}

\bibitem{Andrews80} George E. Andrews, \emph{Totally symmetric plane
 partitions}, Abstracts Amer. Math. Soc. {\bf 1} (1980), 415.

\bibitem{APS} George E. Andrews, Peter Paule, and Carsten Schneider,
\emph{Plane Partitions VI. Stembridge's TSPP theorem},
Adv. Appl. Math. {\bf 34} (2005), 709--739.

\bibitem{B} David M. Bressoud, \emph{Proofs and Confirmations},
Math. Assoc. America and Cambridge University Press (1999).

\bibitem{C} Fr\'ed\'eric Chyzak, \emph{An extension of Zeilberger's
fast algorithm to general holonomic functions}, Discrete Mathematics
{\bf 217} (2000), 115--134.

\bibitem{CS} Fr\'ed\'eric Chyzak and Bruno Salvy,
\emph{Non-commutative elimination in Ore algebras proves
multivariate identities}, J. of Symbolic Computation
{\bf 26} (1998), 187--227.

\bibitem{Kauers} Manuel Kauers, Christoph Koutschan, Doron Zeilberger,
 \emph{A Proof of George Andrews' and Dave Robbins' $q$-TSPP Conjecture (modulo
 a finite amount of routine calculations)}, ArXiv:0808.0571.
 Also available at http://www.math.rutgers.edu/\~{}zeilberg/pj.html

\bibitem{K} Christian Krattenthaler, \emph{Advanced Determinant Calculus},
S\'em. Lothar. Comb. {\bf 42} (1999), B42q.
The Andrews Festschrift, D. Foata and G.-N. Han (eds.)

\bibitem{Koutschan09} Christoph Koutschan, \emph{Advanced Applications of
the Holonomic Systems Approach}, Ph.D. Thesis, RISC, Johannes Kepler University Linz, Austria, 2009.

\bibitem{Koutschan10a} Christoph Koutschan, \emph{Eliminating Human
Insight: An Algorithmic Proof of Stembridge's TSPP Theorem}.
In \emph{Gems in Experimental Mathematics},
Contemporary Mathematics 517, pp. 219-230. 2010. 

\bibitem{O} Soichi Okada, \emph{On the generating functions for certain
classes of plane partitions}, J. Comb. Theory, Series A
{\bf 53} (1989), 1--23.

\bibitem{Sta1} Richard Stanley, \emph{A baker's dozen of conjectures concerning
plane partitions}. In \emph{Combinatoire \'enum\'erative}, ed.
G. Labelle and P. Leroux. Lecture Notes in Mathematics
{\bf 1234}, 285--293. Springer Verlag, New York.

\bibitem{Sta2} Richard Stanley, \emph{Symmetries of plane partitions}, J. Comb. Theory, Series~A
{\bf 43} (1986), 103--113.

\bibitem{Ste} John Stembridge, \emph{The enumeration of totally
symmetric plane partitions}, Advances in Mathematics {\bf 111},
227--243 (1995).

\bibitem{T} Nobuki Takayama, \emph{Gr\"obner basis, integration and
transcendental functions}. In Proceedings of ISSAC'90 (1990), 152--156.

\bibitem{Z1} Doron Zeilberger,
\emph{A Holonomic systems approach to special functions
identities}, J. of Computational and Applied Math. {\bf 32},
321--368 (1990).

\bibitem{Z3} Doron Zeilberger, \emph{The HOLONOMIC ANSATZ II.
Automatic DISCOVERY(!) and PROOF(!!) of Holonomic Determinant Evaluations},
Annals of Combinatorics {\bf 11} (2007), 241--247.

\end{thebibliography}
%% .tex document. To conform to PNAS requirements, copy the reference listings
%% from your .bbl file and add them to the article .tex file, using the
%% bibliography environment described above.  

%%  Contact pnas@nas.edu if you need assistance with your
%%  bibliography.

% Sample bibliography item in PNAS format:
%% \bibitem{in-text reference} comma-separated author names up to 5,
%% for more than 5 authors use first author last name et al. (year published)
%% article title  {\it Journal Name} volume #: start page-end page.
%% ie,
% \bibitem{Neuhaus} Neuhaus J-M, Sitcher L, Meins F, Jr, Boller T (1991) 
% A short C-terminal sequence is necessary and sufficient for the
% targeting of chitinases to the plant vacuole. 
% {\it Proc Natl Acad Sci USA} 88:10362-10366.

%% Enter the largest bibliography number in the facing curly brackets
%% following \begin{thebibliography}

\end{article}

\setcounter{page}{2199}

%%%%%%%%%%%%%%%%%%%%%%%%%%%%%%%%%%%%%%%%%%%%%%%%%%%%%%%%%%%%%%%%

%% Adding Figure and Table References
%% Be sure to add figures and tables after \end{article}
%% and before \end{document}

%% For figures, put the caption below the illustration.
%%
%% \begin{figure}
%% \caption{Almost Sharp Front}\label{afoto}
%% \end{figure}

%% For Tables, put caption above table
%%
%% Table caption should start with a capital letter, continue with lower case
%% and not have a period at the end
%% Using @{\vrule height ?? depth ?? width0pt} in the tabular preamble will
%% keep that much space between every line in the table.

%% \begin{table}
%% \caption{Repeat length of longer allele by age of onset class}
%% \begin{tabular}{@{\vrule height 10.5pt depth4pt  width0pt}lrcccc}
%% table text
%% \end{tabular}
%% \end{table}

%% For two column figures and tables, use the following:

%% \begin{figure*}
%% \caption{Almost Sharp Front}\label{afoto}
%% \end{figure*}

%% \begin{table*}
%% \caption{Repeat length of longer allele by age of onset class}
%% \begin{tabular}{ccc}
%% table text
%% \end{tabular}
%% \end{table*}

\end{document}